\documentstyle[12pt,amscd,amssymb,righttag]{amsart} 
\newtheorem{thm}{Theorem}[section]
\newtheorem{corol}[thm]{Corollary}
\newtheorem{lemma}[thm]{Lemma}
\newtheorem{prop}[thm]{Proposition}
\theoremstyle{definition}
\newtheorem{defin}[thm]{Definition}

\theoremstyle{remark}
\newtheorem{remark}[thm]{Remark}
\numberwithin{equation}{section}
\def\rest#1,#2{#1_{\vert #2}}

\def\Pic{\o{Pic}}
\let\o=\operatorname
\def\cpuno{{\Bbb{P}}^1}
\def\iso{\kern.35em{\raise3pt\hbox{$\sim$}\kern-1.1em\to}\kern.3em} 
\def\oskip{\par\vbox to4mm{}\par}
\def\fp{\times_{\cpuno}}
\def\rest#1,#2{#1_{\vert #2}}
\def\F{{\cal F}}
\def\O{{\cal O}}
\def\Pc{{\cal P}}
\def\M{{\cal M}} 
\def\L{{\cal L}}

\def\Z{{\Bbb Z}}
\def\C{{\Bbb C}}
\def\Q{{\Bbb Q}}
\def\E{{\cal E}}
\def\N{{\cal N}}
\def\I{{\cal I}}
\def\U{{\cal U}}
\def\Qc{{\cal Q}}

\def\td{\o{td}}
\def\ch{\o{ch}}
\def\stdp{\sqrt{\td p}}
\def\stdhp{\sqrt{\td \hat p}}
\def\pic#1{\o{Pic}(#1)}
\def\bT{{\bold T}}\def\bS{{\bold S}}
\def\picfun{{\bold Pic}^-_{X/\cpuno}}
\def\jac{{J(\rest X,{\U})}}
\def\bysame{$\relbar\joinrel\joinrel\relbar\joinrel\joinrel\relbar \joinrel 
\joinrel\relbar\joinrel\joinrel \relbar\joinrel\joinrel
\relbar\joinrel\joinrel\relbar\joinrel\joinrel \relbar\joinrel\joinrel
\relbar$, \ }
\def\X{{\widehat X}}

\def\f{{\bold f}}
\begin{document}
\pagestyle{plain}
\title[]{MIRROR SIMMETRY ON K3 SURFACES \\ VIA FOURIER-MUKAI TRANSFORM} 
\author[]{Claudio Bartocci\thinspace \ddag \\ Ugo Bruzzo\thinspace \S \\ Daniel
Hern\'andez Ruip\'erez \thinspace \P \\ Jos\'e M.~Mu\~noz Porras\thinspace \P}
\address{\ddag\thinspace Dipartimento di Matematica, Universit\`a di
Ge\-no\-va, Via Dodecaneso 35, 16146 Ge\-no\-va, Italy}
\email{bartocci@@dima.unige.it}
\address{\S\thinspace Scuola Internazionale Superiore di Studi Avanzati \\  Via
Beirut 2-4, 34014 Trieste, Italy} \email{bruzzo@@sissa.it}
\address{\P\thinspace Departamento de Matem\'atica Pura y Aplicada,
Universidad de Salamanca, Plaza de la Merced 1-4, 37008 Salamanca, Spain}
\email{ruiperez@@gugu.usal.es, jmp@@gugu.usal.es} 
\date{29 April 1997. Revised 17 October 1997} \subjclass{14F05, 14J27, 14J28,
83E30}
\keywords{Fourier-Mukai transform, K3 surfaces, mirror symmetry} 
\rightline{SISSA Ref. 61/97/fm/geo}
\rightline{alg-geom/9704023}
\rightline{To appear in {\it Commun.~Math.~Phys.}}
\bigskip\bigskip\bigskip\bigskip
\begin{abstract} We use a relative Fourier-Mukai transform on elliptic K3 
surfaces $X$ to describe mirror symmetry. The action of this Fourier-Mukai
transform on the cohomology ring of $X$ reproduces relative T-duality and
provides an infinitesimal isometry of the moduli space of algebraic structures
on $X$ which, in view of the triviality of the quantum cohomology of K3
surfaces, can be interpreted as mirror symmetry.

{From} the mathematical viewpoint the novelty is that we exhibit another 
example of a Fourier-Mukai transform on K3 surfaces, whose properties are
closely related to the geometry of the relative Jacobian of $X$. \end{abstract}
\maketitle
\oskip\section{Introduction} In a recent approach of Strominger, Yau and 
Zaslow \cite{SYZ}, the phenomenon of mirror symmetry on Calabi-Yau threefolds
admitting a
$T^3$ fibration is interpreted as T-duality on the $T^3$ fibres. According  to
this formulation one would like to {\it define} the mirror dual to a Calabi-Yau
manifold (of any dimension) as a compactification of the moduli space of its
special Lagrangian submanifolds (the $T^3$ tori in the above case) endowed with
a suitable complex structure \cite{SYZ,Morr,Hit}.

In two dimensions this means that one considers a K3 surface elliptically
fibred over
the projective line, $p \colon X \to
\cpuno$. A mirror dual to $X$ can be identified with the component $\M$ of
the moduli space of simple sheaves on $X$ having Mukai vector $(0,\mu,0)\in
H^\bullet(X,\Z)$, where
$\mu$ is the cohomology class defined by the fibres of $p$. The mirror map
between the Hodge lattices of $X$ and $\M$ should be given by a suitable
Fourier-Mukai transform \cite{Morr,eng,DM}. 

In this paper we show that a Fourier-Mukai transform on elliptically fibred
K3 surfaces provides indeed a description of mirror symmetry. The
Fourier-Mukai transform not only maps special Lagrangian 2-cycles to
0-cycles, as noticed by Morrison and others, but also reproduces the
correct duality transformations on 4-cycles and on 2-cycles of genus 0. It
turns out that the Fourier-Mukai transform does not define an automorphism
of the cohomology ring of the K3 surface which swaps the directions
corresponding to complex structures with the directions corresponding to
complexified K\"ahler structures. In this sense our treatment is different
from other approaches, cf.~e.g.~\cite{Asp,G,GW}. However, we are able to
obtain an isometry between the tangent space to the deformations of complex
structures on $X$ and the tangent space to the deformations of
``complexified K\"ahler structures'' on the mirror manifold.

We also note that the map determined by the Fourier-Mukai transform has a
correct action on the mass of the so-called BPS states. 

In order to describe this ``geometric mirror symmetry'' two modifications
must be introduced in the construction we have above outlined. First, we
regard the mirror dual to the elliptic K3 surface $X$ as its compactified
relative Jacobian $\X$ (this is actually isomorphic to $\M$); secondly, we
define a Fourier-Mukai in a relative setting (cf.~\cite{Muk3} for a
relative Fourier-Mukai transform for abelian schemes). Moreover, the
relative transform we define, once restricted to the smooth fibres, reduces
to the usual Fourier-Mukai transform for abelian varietes; in this way the
reduction of mirror symmetry to relative T-duality in the spirit of
\cite{SYZ} is achieved. 

It should be stressed that this analysis shows that the moduli space $\M$
is isomorphic to the original K3 surface $X$ as an algebraic variety, in
accordance with the fact that, under this interpretation of mirror
symmetry, a K3 surface is mirror to itself \cite{SYZ}. This fact, together
with the the existence of an isometry between the above mentioned spaces of
deformations, is consistent with the triviality of the quantum cohomology
of a K3 surface (in particular, the Weil-Petersson metric on the moduli
space of complexified K\"ahler structures bears no instantonic
corrections). 

To go through some more detail, the Fourier-Mukai functor $\bT $ we define
transforms a torsion-free rank-one zero-degree sheaf concentrated on an
elliptic fibre of $X$ to a point of the compactified relative Jacobian
$\X$; accordingly, $\bT $ enjoys the T-duality property of relating
$2$-cycles to $0$-cycles. Furthermore, $\bT $ induces an isometry
$$\psi\colon H^{1,1}(\X,\C)/\o{Pic}(\X)\otimes\C\to H^{1,1}( X,\C)/ \o{Pic}( X)
\otimes\C\,.$$ The quotient $H^{1,1}(\X,\C)/\o{Pic}(\X)\otimes\C$ can be
regarded as the tangent space at
$\X$ to the space of deformations of algebraic structures on $\X$ which
preserve the Picard lattice, and similarly, $H^{1,1}( X,\C)/\o{Pic}( X)
\otimes\C$ is to be identified with the tangent space to the space of
deformations of K\"ahler structures on $X$ preserving the Picard lattice.
With these identifications in mind, the isometry $\psi$ can be regarded as
an ``infinitesimal'' mirror map. 

{From} a mathematical viewpoint the transform we define here provides
another example of a 
Fourier-Mukai transform on K3 surfaces in addition to the one given in
\cite{BBH}.

The paper is organized as follows. In Section 2 we fix notations, define
the relative 
Fourier-Mukai functor and prove its first properties. In Section 3 we prove
that it is
invertible and thus gives rise to an equivalence of derived categories. In
Section 4 we study
the action of the Fourier-Mukai transform on the cohomology ring of the K3
surface $X$. In
Section 5 we discuss how the Fourier-Mukai transform can be regarded as a
mirror duality for
string theories compactified on an elliptic K3 surface.

\oskip\section{The basic construction}
Let $p\colon X\to \cpuno$ be a minimal algebraic elliptically fibred K3
surface (all 
algebraic varieties will be over $\C$). The fibration $p$ has singular
fibres; these has been
classified by Kodaira \cite{Kod}. We assume that
$p\colon X\to\cpuno$ has a section $e\colon\cpuno\hookrightarrow X$ and
write $H=e(\cpuno)$. 
We shall denote by $X_t$ the fibre of $p$ over $t\in\cpuno$, and by $i_t\colon
X_t\hookrightarrow X$ the inclusion. 

\smallskip
{\it A compactification of the relative Jacobian.} Let $\M$ be the moduli 
space of simple sheaves on $X$, of pure dimension 1 and Chern character
$(0,\mu,0)$, where $\mu$ is the cohomology class of the fibres of $p$. Results
of Simpson \cite{Simp} imply that 
$\M$ is a smooth projective surface  (actually, a minimal  K3 surface,
cf.~\cite{Muk3}).
One may define a morphism
\begin{align}
\gamma\colon X & \to \M \\ x & \mapsto (i_t)_\ast ({\frak m}_x\otimes
\O_{X_t}(e(t)))
\label{e:gamma}\end{align} where $X_t\ni x$, and ${\frak m}_x$ is the ideal
sheaf of $x$  in
$X_t$.

Let $\U\subset\cpuno$ be the open subset supporting the smooth fibres of
$p$, and let
$\jac$ be the relative Jacobian variety. The restriction of $\gamma$ to
$\rest X,{\U}$ 
factors as
$$\rest X,{\U} \iso \jac\hookrightarrow \M\,, $$ where the isomorphism
$\rest X,{\U}  \iso
\jac$ is given by $x\mapsto \O_{X_t}(e(t)-x)={\frak
m}_x\otimes\O_{X_t}(e(t))$ and
$\jac\hookrightarrow
\M$ associates with any zero-degree torsion-free sheaf $L_t$ over $X_t$ its
direct image 
$(i_t)_\ast L_t$. Then $\gamma$ is birational, and $X\simeq\M$ since they
both are smooth
projective surfaces and $X$ is minimal.

We want now to construct a suitable compactification of the relative
Jacobian of $p\colon 
X\to \cpuno$. We denote by $\picfun$ the functor which to any morphism
$f\colon S\to\cpuno$ of
algebraic varieties associates the space of $S$-flat sheaves on $p_S\colon
X\fp S\to S$, whose
restrictions to the fibres of $p_S$ are torsion-free, of rank one and
degree zero. Two such
sheaves
$\F$, $\F'$ are considered to be equivalent if $\F'\simeq\F\otimes
p_S^\ast\L$ for a line 
bundle
$\L$ on $S$ (cf.~\cite{AK}). Due to the
existence of the section $e$, $\picfun$ is a sheaf functor. 

\begin{prop} The functor $\picfun$ is represented by an algebraic variety
$\hat p\colon\X
\to\cpuno$, which is isomorphic to $X$. \end{prop} \begin{pf} If we denote
by ${\bold h}_X$, ${\bold h}_{\M}$ the functors of points of $X$, $\M$ as
schemes over $\cpuno$, the isomorphism $\gamma\colon {\bold h}_X\iso {\bold
h}_{\M}$ factors as
$$ {\bold h}_X@>\varpi>>\picfun @>\alpha>> {\bold h}_{\M} $$ where $\varpi$
and $\alpha$ are defined (over the closed points) by $\varpi(x)= {\frak
m}_x\otimes\O_{X_t}(e(t))$ and $\alpha(L_t)=(i_t)_\ast L_t$ for any
zero-degree torsion-free sheaf $L_t$ over $X_t$. Both morphisms of functors
are immersions and their composition is an isomorphism, so that they are
isomorphisms as well. Then, $\picfun$ is represented by a fibred algebraic
variety $\hat p\colon\X\to\cpuno$, and $\varpi\colon X\to\X$ is an
isomorphism. \end{pf} We denote by $\hat e\colon\cpuno\to\X$ the canonical
section; one has
$\hat e=\varpi\circ e$. Moreover,
we denote by $\pi$, $\hat\pi$ the projections of the fibred product $X\fp
\X$ onto the factors.

\begin{remark} The Picard functor is also representable by an open dense
subscheme
$J$ of $\X$, the relative Jacobian $J\to\cpuno$ of $X\to\cpuno$. If
$X^s\subset X$ denotes the 
complement of the singular points of the  fibres of $\pi$, then $\varpi$
gives an isomorphism
$X^s\iso J$ of schemes over $\cpuno$. One should notice that in general the
Jacobian variety
$J\to\cpuno$ is different from
$\o{Pic}^0(X/\cpuno)\to\cpuno$. This scheme can be obtained from
$J\to\cpuno$ by removing the
images by $\varpi$ of those components of the singular fibres of $p\colon
X\to\cpuno$ that do
not meet the image $H$ of the section. \qed\end{remark}

The representability of $\picfun$ means there exists a coherent sheaf $\Pc$
on $X\fp\X$ flat over
$\X$, whose restrictions to the fibres of $\hat\pi$ are torsion-free, and
of rank one and degree zero, such that \begin{equation}
\o{Hom}\,_{\cpuno}(S,\X) \to \picfun(S)\,,\qquad\qquad \phi \mapsto
[(1\times\phi)^\ast\Pc]
\label{e:univ}\end{equation} is an isomorphism of functors. $\Pc$ is
defined up to tensor product by the pullback of a line bundle on $\X$, and
is called the universal Poincar\'e sheaf. 

To normalize the Poincar\'e sheaf we notice that the isomorphism
$\varpi\colon X\iso\X$ is induced, according to the universal property
(\ref{e:univ}), by the sheaf $\I_\Delta\otimes p_1^\ast\O_X(H)$ on $X\fp
X$, where $p_1$ is the projection onto the first factor and $\I_\Delta$ is
the ideal sheaf of the diagonal $\delta\colon X\hookrightarrow X\fp X$;
this sheaf is flat over the second factor and has zero relative degree.
Then \begin{equation}
\Pc=(1\times\varpi^{-1})^\ast\left(\I_\Delta\otimes p_1^\ast\O_X(H)\right)
\otimes\hat \pi^\ast\L
\label{e:poin}
\end{equation}
for a line bundle $\L$ on $\X$. Restriction to $H\fp\X$ gives
$\rest{\Pc},{H\fp\X}\simeq
\O_{\X}(-2)\otimes\L$, which is a line bundle. We can then normalize $\Pc$
by letting
\begin{equation}
\rest\Pc,{H\fp \X}\simeq \O_{\X}\,.\label{e:norm}\end{equation} We shall
henceforth assume that $\Pc$ is normalized in this way. We shall denote by
$\Pc_\xi$ the restriction
$\rest{\Pc},{\hat\pi^{-1}(\xi)}$. As a consequence of (\ref{e:poin}), $\Pc$
is flat over
$X$ as well.

\begin{remark} Since the moduli space $\M$ is fine, on $X\times \M\simeq
X\times \X$ there is 
a universal sheaf $\Qc$. This is the sheaf that gives rise to the morphism
$\gamma$ (cf.~Eq.~(\ref{e:gamma})). One can show that $\Qc$ is supported on
$X\fp\X$, and its 
restriction to its support is isomorphic to $\Pc$ (up to tensoring by a
pullback of a line
bundle on $\X$).\qed\end{remark}

The dual $\Pc^\ast$ of the Poincar\'e bundle is a coherent sheaf on
$X\fp\X$ whose restrictions to the 
fibres of $\hat\pi$ are torsion-free, rank one, and of degree zero. As
$\Pc^\ast$ is flat over
$\X$ it defines a morphism $\iota\colon\X\to \X$. Since
${\cal E}xt\, ^1_{\O_{X_t}}(\Pc_\xi,\O_{X_t})=0$ for every point $\xi\in\X$
(here $t=\hat 
p(\xi)$) and ${\cal E}xt\,^1_{\O_{X\fp\X}}(\Pc,\O_{X\fp\X})=0$ by
(\ref{e:poin}), the base
change property for the local Ext's (\cite{AK}, Theorem 1.9) implies that
$(\Pc^\ast)_\xi\simeq (\Pc_\xi)^\ast$. Then, the morphism $\iota
\colon\X\to \X$ maps  any
rank-one torsion-free and zero-degree coherent sheaf $\F$ on a fibre $X_t$
to its dual
$\F^\ast$. By (\ref{e:univ}) one has
$(1\times\iota)^\ast\Pc\simeq\Pc^\ast\otimes\hat\pi^\ast\N$ for some line
bundle $\N$ on
$\X$, which turns out to be trivial by (\ref{e:norm}). Then \begin{equation} 
(1\times\iota)^\ast\Pc\simeq\Pc^\ast\,. \label{e:inv}\end{equation} The morphism
$\iota\circ\iota\colon\X\to\X$ is the identity on the Jacobian
$\jac\subset\X$; by
separateness $\iota\circ\iota=\o{Id}$, and (\ref{e:inv}) implies
$\Pc\simeq\Pc^{\ast\ast}\,.$
Then, every coherent sheaf $\F$ on $X\fp S$ flat over $S$ whose
restrictions to the  fibres of
$X\fp S \to S$ are torsion-free and of rank one and degree zero is reflexive,
$\F\simeq\F^{\ast\ast}$.

\begin{prop} The relative Jacobian of the fibration $\hat p\colon\X\to
\cpuno$ admits a 
compactification which is isomorphic to $X$ as a fibred variety over
$\cpuno$, and the
relevant universal Poincar\'e sheaf may be identified with $\Pc^\ast$.
\end{prop} \begin{pf}
By (\ref{e:poin}), the sheaf $\Pc^\ast$ is flat over $X$. Proceeding as above, one proves that
$(\Pc^\ast)_x\simeq (\Pc_x)^\ast$ for every point $x\in X$, which means
that the restrictions
of $\Pc^\ast$ to the fibres of $\pi$ are torsion-free sheaves of rank one
and degree zero. So
$\Pc^\ast$ defines a morphism $X\to \widehat{\X}$ of schemes over
$\cpuno$. If
$\U\subset\cpuno$ denotes as above the open subset supporting the smooth
fibres of  $p$, this
morphism restricts to an isomorphism $X_{\vert\U}\simeq J(\X)_{\vert\U}$ of
schemes over
${\U}$. Since $\X$ is minimal and $X$ is smooth and has no $(-1)$-curves,
$X\simeq\widehat{\X}$.
\end{pf}

The roles of $X$ and
$\X$ are so completely interchangeable.

\par\smallskip
{\it The Fourier-Mukai functors.} For any morphism $f\colon S\to\cpuno$ let
us consider the 
diagram
$$\begin{CD} (X\fp\X)_S @>\hat\pi_S>> \X_S \\ @V\pi_SVV @VV\hat p_SV \\ X_S
@>p_S>> S
\end{CD}
$$ We shall systematically denote objects obtained by base change to $S$ by
a subscript $S$. 
(Note that $(X\fp\X)_S\simeq X_S\times_S\X_S$). We define the Fourier-Mukai
functors $\bS
_S^i$, $i=0,1$ by associating with every sheaf $\F$ on $X_S$ flat over $S$
the sheaf on $\X_S$
$$\bS_S^i(\F)=R^i\hat\pi_{S\ast}(\pi_S^\ast\F\otimes\Pc_S)\,.$$ 

The Fourier-Mukai functors mapping sheaves on $X$ to sheaves on $\X$ will
be denoted by 
$\bS^i$.

\begin{defin} We say that a coherent sheaf $\F$ on $X_S$ flat over $S$ is
WIT$_i$ if
$\bS _S^j(\F)=0$ for $j\ne i$. We say that $\F$ is IT$_i$ if it is WIT$_i$
and  $\bS _S^i(\F)$
is locally free. \end{defin} One should notice that, due to the presence of
the fibred instead
of the cartesian product, the WIT$_0$ and IT$_0$ conditions are not
equivalent: for instance
$\kappa(x)$ (the skyscraper sheaf concentrated at $x\in X$) is WIT$_0$ but
not $IT_0$.

Since the fibres of $\hat\pi_S$ are one-dimensional the first direct image
functor commutes with base change. \begin{prop} Let $\F$ be a sheaf on
$X_S$, flat over $S$. For every morphism $g\colon T\to S$ one has
$g_{\X}^\ast\bS_S^1(\F)\simeq \bS_T^1(g_X^\ast\F)$, where
$g_X\colon X_T\to X_S$, $g_{\X}\colon\X_T\to\X_S$ are the morphisms induced
by $g$.
\qed\label{basechange}\end{prop}

The zeroth direct image does not commute with base change; however, a
weaker property holds.

\begin{prop} Let $\F$ be a sheaf on $X_S$, flat over $S$. For every point
$\xi\in \X_S$, the natural base change morphism $$
\hat\pi_{S\ast}(\pi_S^\ast\F\otimes\Pc_S)\otimes\kappa(\xi)\to
H^0(X_s,\F_s\otimes\Pc_\xi)
$$
is injective (here $s=\hat p_S(\xi)$\rm). \label{basechange2}
\end{prop}
\begin{pf} Let ${\frak m}_\xi$ denote the ideal sheaf of $\xi\in \X_S$.
Since $\hat\pi_S$ is flat, $\hat\pi_S^\ast{\frak m}_\xi $ is the ideal
sheaf of the fibre $\hat\pi_S^{-1}(\xi)\simeq X_s$ in $X_S\times_S\X_S$.
Let us write $\N=\pi_S^\ast\F\otimes\Pc_S$ and
$\N_\xi=\rest{\N},{\hat\pi_S^{-1}(\xi)}$. Since $\N$ is flat over $\X_S$
there is an exact sequence $$ 0@>>>\hat\pi_S^\ast{\frak m}_\xi \otimes\N
@>>> \N @>>> j_\ast \N_\xi @>>> 0\,, $$ where $j\colon
\hat\pi_S^{-1}(\xi)\simeq X_s\hookrightarrow X_S\times_S\X_S$ is the
natural immersion. By taking direct images we obtain $$ 0 @>>>
\hat\pi_{S\ast} (\hat\pi_S^\ast{\frak m}_\xi \otimes\N) @>>>
\hat\pi_{S\ast}(\N) @>\eta>>
\hat\pi_{S\ast}(j_\ast\N_\xi)=H^0(X_s,\N_\xi)\,. $$ By the projection
formula, $\hat\pi_{S\ast} (\hat\pi_S^\ast{\frak m}_\xi \otimes\N)\simeq
{\frak m}_\xi\otimes \hat\pi_{S\ast}\N$, and then
$\ker\eta\simeq {\frak m}_\xi\cdot\N$; this implies that the base change
morphism
$\N\otimes\kappa(\xi)\to H^0(X_s,\N_\xi)$ is injective.\end{pf} 

\smallskip {\it Fourier-Mukai transform of rank 1 sheaves.} A first
manifestation of geometric mirror symmetry is the fact that the
Fourier-Mukai transform of a torsion-free rank-one zero-degree coherent
sheaf on a fibre $X_t$ is a skyscraper sheaf concentrated at a point of
$\X_t$.

By Proposition \ref{basechange} the basic ingredients to compute the
functors $\bS^\bullet$ are the Fourier-Mukai transforms
$\bS_{\X}^\bullet(\Pc)$ of the universal Poincar\'e sheaf $\Pc$ on
$X_{\X}=X\fp\X$. The relevant higher direct images of $\Pc$ and $\Pc^\ast$
are computed as follows. (For every algebraic variety $q\colon Y\to\cpuno$
over $\cpuno$ and every coherent sheaf $\N$ on $Y$ we denote by $\N(n)$ the
sheaf $\N\otimes q^\ast\O_{\cpuno}(n)$.) \begin{thm} $\phantom{xxxxx}$
\begin{list}{}{\itemsep=2pt}
\item[1] $\bS_{\X}^1(\Pc)\simeq
\zeta_\ast\O_{\X}(-2)$, where $\zeta\colon \X\hookrightarrow \X\fp\X$ is
the graph of the morphism $\iota$.
\item[2.] $\bS_{\X}^0(\Pc)=0$.
\item[3.] $\bS_{\X}^1(\Pc^\ast)\simeq \delta_\ast\O_{\X}(-2),\quad
\bS_{\X}^0(\Pc^\ast)=0$, where $\delta\colon\X\hookrightarrow\X\fp\X$ is
the diagonal immersion.
\item[4.] $R^1\hat\pi_\ast\Pc\simeq R^1\hat\pi_\ast\Pc^\ast \simeq \hat
e_\ast\O_{\cpuno}(-2)$, while the zeroth direct images vanish. \end{list}
\label{pm}\end{thm}
A result similar to the second formula can be found in \cite{Muk3} for the
case of relative abelian schemes.

To prove Theorem \ref{pm} we need some preliminary results. 

\begin{lemma} Let $Y$ be a proper connected curve of arithmetic genus 1 and
$\F$ a torsion-free rank-one and zero-degree sheaf on $Y$. Then
$H^1(Y,\F)\neq 0$ if and only if $\F\simeq\O_Y$. \label{huno}
\end{lemma}
\begin{pf} One has $H^0(Y,\F)\neq 0$ by Riemann-Roch and $H^0(Y,\F^\ast)
\neq 0$ by duality. Let $\tau$ and $\sigma$ be nonzero sections of $\F$ and
$\F^\ast$ respectively. Let
$\rho$ be the composition $$\begin{CD} \F @>>> \F^{\ast\ast} @>
\sigma^\ast>> \O_Y\,.\end{CD}$$ Since
$\rho\circ\tau\ne 0$, the morphism $\rho\circ\tau$ consists in the
multiplication by a nonzero constant, which may be set to 1. Then
$\rho\circ\tau=\o{id}$, so that
$\F\simeq
\O_Y\oplus\M$, where $\M$ has rank zero; hence $\M=0$, and $\F\simeq\O_Y$.
\end{pf}

\begin{lemma} Let $\xi,\mu\in \X_t$. \begin{list}{}{\itemsep=2pt} 
\item[1.]  
The sheaf
$\Pc_\xi\otimes\Pc_\mu$ has torsion if and only if $\xi$ is a singular
point of the fibre
$\X_t$ and $\mu=\xi$.

\item[2.] The evaluation morphism $\Pc_\xi\otimes\Pc_\xi^\ast\to\O_{X_t}$
induces
an isomorphism
$H^1(X_t,\Pc_\xi\otimes\Pc_\xi^\ast)\simeq H^1(X_t,\O_{X_t})$. 

\item[3.] If $\mu\ne\iota(\xi)$ then $H^1(X_t,\Pc_\xi\otimes\Pc_\mu)=0$.
\end{list}
\label{eval}
\end{lemma}
\begin{pf}
1. We have $\Pc_\xi={\frak m}_x(e(t))$ for a point $x\in X_t$, and $\xi$ 
is singular in $\X_t$ iff $x$ is singular in $X_t$. 
If $\xi$ or $\mu$ are not singular, then
one of the sheaves $\Pc_\xi$, $\Pc_\mu$ is locally free, and
$\Pc_\xi\otimes\Pc_\mu$ is
torsion-free. 
Otherwise, 
$\Pc_\xi={\frak m}_x(e(t))$ and $\Pc_\mu={\frak m}_y(e(t))$ for singular points
$x,y\in X_t$. If
$\mu\ne\iota(\xi)$ then
${\frak m}_x\otimes{\frak m}_y$ is torsion-free. 
Finally, if $\mu=\iota(\xi)$ then
${\frak m}_x(2e(t))={\frak m}_y^\ast$, so that $x=y$ (because ${\frak
m}_x$, ${\frak m}_y$ are
not locally-free only at $x$, $y$, respectively). Thus $\mu=\xi$, and
${\frak m}_x\otimes{\frak
m}^\ast_x$ has torsion at $x$. 

2. The only nontrivial case is when $\Pc_\xi$ is not locally free. 
Let $x\in X_t$ be the singular point corresponding to $\xi$. We have an
exact sequence
$$
0@>>>{\frak m}_x@>>>\Pc_\xi\otimes\Pc_\xi^\ast @>>>{\frak m}_x/{\frak
m}_x^2 @>>>0\,,
$$
which implies that $H^1(X_t,{\frak m}_x)\to
H^1(X_t,\Pc_\xi\otimes\Pc_\xi^\ast)$ is an epimorphism. Since the
composition $H^1(X_t,{\frak m}_x)\to H^1(X_t,\Pc_\xi\otimes\Pc_\xi^\ast)
\to H^1(X_t,\O_{X_t})$ is an isomorphism,
$H^1(X_t,\Pc_\xi\otimes\Pc_\xi^\ast)\simeq H^1(X_t,\O_{X_t})$ is an
isomorphism as well.

3. Follows from 1 and Lemma \ref{huno}.
\end{pf}

In order to compute the Fourier-Mukai transform $\bS_{\X}^\bullet(\Pc)$ of
the Poincar\'e sheaf
$\Pc$ on $X\fp \X$ we consider the diagram $$ \begin{CD} X\fp\X\fp\X
@>\pi_{23}>> \X\fp\X \\ @V\pi_{12}VV @VV\hat p_1V \\ X\fp\X @>\hat \pi >>
\X
\end{CD}
$$ The Poincar\'e sheaf on $X\fp\X\fp\X$ is $\Pc_{\X}=\pi_{13}^\ast\Pc$ and
the Fourier-Mukai transforms of $\Pc$ are
$\bS_{\X}^\bullet(\Pc)=R^\bullet\pi_{23\ast}(\pi_{12}^\ast\Pc\otimes
\pi_{13}^\ast\Pc)$.

\smallskip{\it Proof of 1 of Theorem \ref{pm}.} We have $\Pc\otimes\Pc^\ast=
(1\times\zeta)^\ast(\pi_{12}^\ast\Pc\otimes\pi_{13}^\ast\Pc)$. The composition of the epimorphism
$\pi_{12}^\ast\Pc\otimes\pi_{13}^\ast\Pc@>>>(1\times\zeta)_\ast
(\Pc\otimes\Pc^\ast)$ with the evaluation morphism
$(1\times\zeta)_\ast(\Pc\otimes\Pc^\ast)@>>>(1\times\zeta)_\ast
(\O_{X\fp\X})$ gives a morphism
$$
\pi_{12}^\ast\Pc\otimes\pi_{13}^\ast\Pc@>>>(1\times\zeta)_\ast (\O_{X\fp\X})\,.
$$ We have then a morphism
$$
\bS_{\X}^1(\Pc)=R^1\pi_{23\ast}(\pi_{12}^\ast\Pc\otimes\pi_{13}^\ast\Pc)
@>>> R^1\pi_{23\ast} ((1\times\zeta)_\ast(\O_{X\fp\X}))\simeq
\zeta_\ast\O_{\X}(-2)\,. $$ Since the first direct image functor commutes
with base change, Lemma \ref{eval} implies that $\bS_{\X}^1(\Pc)$ is
supported on $\zeta(\X)$. Moreover the fibre of the previous morphism at a
point $\zeta(\xi)$, with $\xi\in \X_t$, is
$H^1(X_t,\Pc_\xi\otimes\Pc_\xi^\ast)@>>>H^1(X_t,\O_{X_t})$, which is an
isomorphism by Lemma \ref{eval}.\qed\smallskip 

Let $f\colon S\to\cpuno$ be a morphism and $\F$ a coherent sheaf on $X\fp
S$ flat over
$S$ whose restrictions to the fibres of $p_S$ are torsion-free and have
rank one and degree zero. Let
$\phi\colon S\to\X$ be the morphism determined by the universal property
(\ref{e:univ}), so
that
$$
(1\times\phi)^\ast\Pc\simeq\F\otimes p_S^\ast\L\,, $$ for a line bundle
$\L$ on $S$. Let
$\Gamma\colon S\hookrightarrow \X_S$ be the graph of the morphism
$\iota\circ\phi\colon S\to\X$.

\begin{lemma}
$\bS_S^1(\F)\otimes \hat p_S^\ast\L\simeq \Gamma_{\ast}\O_S(-2)\,,\quad
\bS_S^0(\F)=0$.
\label{SS}\end{lemma}
\begin{pf} The formula for $\bS_S^1(\F)$ follows from Proposition
\ref{basechange} and 1 of Theorem \ref{pm} after some standard
computations. The second formula is proved as follows. From Proposition
\ref{basechange2} we have the exact sequence $$ 0 @>>>
\hat\pi_{S\ast}(\pi_S^\ast\F\otimes\Pc_S)\otimes\kappa(\xi) @>>>
H^0(X_s,\F_s\otimes\Pc_\xi)$$
where $s=\hat p_S(\xi)$.
If $\xi\notin\Gamma(S)$, $H^0(X_s,\F_s\otimes\Pc_\xi)=0$ by Lemma
\ref{eval} and
$\hat\pi_{S\ast}(\pi_S^\ast\F\otimes\Pc_S)\otimes\kappa(\xi)=0$ as well. 
If
$\xi\in\Gamma(S)$ the first direct image $\bS_S^1(\F)$ is not locally-free
at $\xi$ since it is concentrated on the image of $\Gamma$, and then the
second arrow is not surjective; but $H^0(X_s,\F_s\otimes\Pc_\xi)$ is
one-dimensional by Lemma \ref{eval}, so that
$\hat\pi_{S\ast}(\pi_S^\ast\F\otimes\Pc_S)\otimes\kappa(\xi)=0$. \end{pf}

\smallskip{\it End of proof of Theorem \ref{pm}.} 2 is proved by applying
Lemma \ref{SS} with $S=\X$ and $\phi$ the identity, while to prove 3 one
chooses
$S=\X$ and $\phi=\iota$. Taking $S=\cpuno$ and $\phi=\hat e$ one proves the
claims of 4 concerning the sheaf $\Pc$. To prove the claims for $\Pc^\ast$
one notices that Lemma \ref{SS} still applies after replacing $\Pc$ by
$\Pc^\ast$.\qed\smallskip 

We can now compute the Fourier-Mukai transform of sheaves on $X$
corresponding to points in $\X$.
\begin{corol}
Let $\F$ be a rank-one, zero-degree, torsion-free coherent sheaf on a fibre
$X_t$. Then
$$
\bS^0_t(\F)=0\,,\qquad \bS ^1_t(\F)=\kappa([\F^\ast])\,, $$ where
$[\F^\ast]$ is the point of
$\X_t$ defined by $\F^\ast$. \qed\label{cor2}\end{corol} 

\oskip\section{Inversion of the Fourier-Mukai transform} The Fourier-Mukai
functor defines a functor $D^-(X)\to D^-(\X)$ given by $\bS(F)=
R\hat\pi_\ast(\pi^\ast F{\overset L \otimes}\Pc)$ (here $D^-(X)$ is the
subcategory of the derived category of coherent
$\O_X$-modules consisting of the complexes bounded from above). To state
the invertibility properties of this functor in a neat way we define a
modified functor $\bT \colon D^-(X)\to D^-(\X)$ by $\bT (F)=\bS(F{\overset
L\otimes}\O_X(1))$. A natural candidate for the inverse of
$\bT $ is the functor
$\widehat{\bT}\colon D^-(\X)\to D^-(X)$ given by $$\widehat{\bT}(G)=
\widehat{\bS} (G{\overset L\otimes}\O_{\X}(1))\qquad
\text{where}\qquad\widehat{\bS}(G') =R\pi_\ast(\hat\pi^\ast G'{\overset
L\otimes}\Pc^\ast)\,.$$ Since the relative dualizing complexes of $\pi$ and
$\hat\pi$ are both isomorphic to $\O_{X\fp\X}(2)[1]$, relative duality
gives:

\begin{prop} For every objects $F$ in $D^-(X)$ and $G$ in $D^-(\X)$ one has
functorial isomorphisms $$
\o{Hom}_{D^-(\X)}(G,\bT (F))\simeq\o{Hom}_{D^-(X)}(\widehat{\bT}(G), F[-1])
$$ $$
\o{Hom}_{D^-(X)}(F,\widehat{\bT}(G))
\simeq\o{Hom}_{D^-(\X)}(\bT (F),G[-1])\,. $$ \qed\end{prop} 

\begin{thm} For every $G\in D^-(\X)$, $F\in D^-(X)$ there are functorial
isomorphisms
$$
\bT (\widehat{\bT}(G))\simeq G[-1]\,,\quad \widehat{\bT}(\bT (F))\simeq
F[-1] $$ in the derived
categories $D^-(\X)$ and $D^-(X)$, respectively. \end{thm} \begin{pf} Let
$\hat\pi_1$ and $\hat\pi_2$ be the projections onto the two factors of
$\X\fp\X$. Then
$ \bT (\widehat{\bT}(G))=R \hat\pi_{2,\ast}(\hat\pi_1^\ast G{\overset
L\otimes}\widetilde{\Pc})
\otimes\O_X(2)$ (see \cite{Muk1,Muk3,Muk4,BBH} for similar statements), with $$
\widetilde{\Pc}=R\pi_{23,\ast}(\pi_{12}^\ast\Pc^\ast\otimes
\pi_{13}^\ast\Pc)\,. $$ By Theorem
\ref{pm} $\widetilde{\Pc}\simeq\delta_\ast(\O_{\X}(-2))[-1]$ in the derived
category, and
$\bT (\widehat{\bT}(G))\simeq G[-1]$.

The second statement follows from the first by interchanging the roles of
$X$ and $\X$.
\end{pf} So $\bT $ establishes an equivalence of triangulated categories.
\begin{corol} Let $\F$
be a WIT$_i$ sheaf on $X$. Then its Fourier-Mukai transform $\bT^i(F)$ is a
WIT$_{1-i}$ sheaf on
$\X$, whose Fourier-Mukai transform
$$\widehat{\bT}^{1-i}(\bT^i(F))=R^{1-i}\pi_\ast(\hat\pi^\ast \bT^i(F)
\otimes\Pc^\ast(1))$$ is isomorphic to $\F$. \label{inv2} \qed\end{corol}

We also have a property of preservation of the Hom groups, which is
sometimes called ``Parseval theorem.''
\begin{prop} There are functorial isomorphisms $$ \o{Hom}_{D^-(\X)} (G,\bar
G) \simeq
\o{Hom}_{D^-(X)}(\widehat{\bS}(G),\widehat \bS(\bar G))\simeq
\o{Hom}_{D^-(X)} (\widehat{\bT}(G),\widehat
\bT(\bar G)) $$ $$
\o{Hom}_{D^-(X)}(F,\bar F) \simeq
\o{Hom}_{D^-(\X)}(\bS(F),\bS(\bar F))\simeq \o{Hom}_{D^-(\X)}(\bT(F),
\bT(\bar F)) $$ for $F$, $\bar F$ in $D^-(X)$ and $G$, $\bar G$ in
$D^-(\X)$.
\qed\end{prop}
\begin{corol} Let $\F$, $\F'$ be coherent sheaves on $X$. If $\F$ is
WIT$_i$ and $\F'$ is WIT$_j$, we have
$$
\o{Ext}^h(\F,\F')\simeq\o{Ext}^{h+i-j}(\bS^i(\F),\bS^j(\F'))
\simeq\o{Ext}^{h+i-j} (\bT^i(\F),\bT^j(\F'))\,. $$ for $h=0,1$. In
particular, if $\F$ is WIT$_i$ there is an isomorphism
$\o{Ext}^h(\F,\F)\simeq\o{Ext}^h(\bT^i(\F),\bT^i(\F))$ for every $h$, so
that
$\bT^i(\F)$ is simple if $\F$ is.
\qed\label{simple}\end{corol}

\begin{remark} Moduli spaces of sheaves on holomorphic symplectic surfaces
carry a holomorphic symplectic structure, which is given by the Yoneda
pairing $\o{Ext}^1(\F,\F)\otimes
\o{Ext}^1(\F,\F)\to\o{Ext}^2(\F,\F)\simeq \C$ (cf.~\cite{Muk2}), where one
identifies
$\o{Ext}^1(\F,\F)$ with the tangent space to the moduli space at the point
corresponding to the
sheaf $\F$. Whenever the Fourier-Mukai transform establishes a morphism
between such moduli
spaces, Corollary \ref{simple} implies that the morphism is symplectic.
\qed\label{sympl}\end{remark}

\oskip\section{Action on the cohomology ring} The cohomology ring
$H^\bullet(X,\Z)$ carries a bilinear pairing, usually called {\it Mukai
pairing,} defined as
$$(a,b,c)\cdot(a',b',c')=(b\cup b'-a\cup c'-a'\cup c)\,\backslash\,[X], $$
and the same is true for $H^\bullet(\X,\Z)$ (here $\,\backslash\,$ denotes
the slant product). We define an isomorphism $\f\colon H^\bullet(X,\Q)\to
H^\bullet(\X,\Q)$ and want to show that in terms of $\f $ one can introduce
an isometry between the tangent space to the moduli space of algebraic
structures on $\X$ and the space of deformations of the complexified
K\"ahler structure on $X$, which can be regarded as a geometric realization
of mirror symmetry. We define the map $\f $ basically as in \cite{Muk4},
but the properties of this map are slightly different, since we are working
in a relative setting, and the relative dualizing sheaf is nontrivial.
Also, we must take coefficients in $\Q$ because the relative Todd
characters involved in the definition of the $\f$ map do not have integral
square roots. 

\smallskip {\it The $\f $ map.} We now define the $\f $ map and describe
its basic properties. We shall be concerned with varieties fibred over
$\cpuno$, $\phi_Y\colon Y\to\cpuno$, with a section $\sigma_Y\colon \cpuno
\hookrightarrow Y$. Since
$\sigma_Y^\ast\circ\phi_Y^\ast=1$, there is a decomposition
$$H^\bullet(Y,\Q)\simeq \phi_Y^\ast H^\bullet(\cpuno,\Q)\oplus
H^\bullet_\phi(Y,\Q)
$$ where $H^\bullet_\phi(Y,\Q)=\ker\sigma_Y^\ast$. One has in particular $$
H^0_\phi(Y,\Q)=0,
\quad H^2 (Y,\Q)=\Q\mu_Y\oplus H^2_\phi(Y,\Q),\quad H^{2i}_\phi(Y,\Q)=
H^{2i}(Y,\Q)\quad \text{for
$i\ge2$}\,. $$ We define in $H^{\text{even}} (Y,\Q)$ an involution $^\ast$
by letting
\begin{equation}\begin{align*} \alpha^\ast =(-1)^i\alpha&\qquad
\text{if}\quad \alpha\in H^{2i}_\phi(Y,\Q) \\ (\phi_Y^\ast\eta)^\ast
=\phi_Y^\ast\eta&\qquad\text{if}\quad \eta\in H^{2i}(\cpuno,\Q)
\,.\end{align*} \end{equation}

Turning back to the case where $X$ is an elliptic K3 surface, satisfying
all the properties we have so far stated, we define morphisms $$\f \colon
H^\bullet(X,\Q)\to H^\bullet(\X,\Q),\qquad \f '\colon H^\bullet(\X,\Q)\to
H^\bullet(X,\Q)$$ by letting $$\f
(\alpha)=\hat\pi_\ast(Z\,\pi^\ast\alpha),\qquad \f
'(\beta)=\pi_\ast(Z^\ast\,\hat\pi^\ast\beta)\,.$$ where
$$Z=\sqrt{\td\hat\pi}\,\ch(\Pc\otimes\pi^\ast\O_X(1))\, \sqrt{\td\pi}\,.$$
\begin{lemma} The maps $\f $, $\f '$ have the following properties:
\begin{list}{}{\itemsep=2pt}
\item[1.] $\f \circ \f '(\beta)=-\beta$; \item[2.] $\f $ and $\f '$ are
$H^\bullet(\cpuno,\Q)$-module isomorphisms; \item[3.] $\f (\mu)=-\hat w$,
where $\hat w$ is the fundamental class of $\X$; \item[4.] $\f (H)=1+\hat
w$;
\item[5.] $\f (1)=-\hat\mu-\Theta+\hat w$, where $\hat \mu$ is the divisor
given by the 
fibres of $\hat p\colon \X\to\cpuno$, and $\Theta=\hat e(\cpuno)$.
\item[6.] $\beta\cdot \f (\alpha)= -\f'(\beta)\cdot\alpha$ for $\alpha\in
H^\bullet_p(X,\Q)$,
$\beta\in H^\bullet_{\hat p}(\X,\Q)$.
\item[7.] $\f $ establishes an isometry between $H^\bullet_p(X,\Q)$ and
$H^\bullet_{\hat p}(\X,\Q)$.
\end{list}
\label{newlemma}\end{lemma}
\begin{pf} Property 1 is proved as in \cite{Muk4}, p.~382, provided that
suitable adaptations to the relative case are done. One also proves that
$\f '\circ \f (\alpha)=-\alpha$, so that 2 follows. 

To prove 3, let $\L$ be a flat line bundle on a smooth fibre $X_t$ of $p$.
One knows that $\ch i_{t\ast}(\L)=i_{t\ast}(1)=\mu$ since the normal bundle
to $X_t$ is trivial. By Corollary \ref{cor2} we have $$\bS
^0(i_{t\ast}\L)=0,\qquad\bS ^1(i_{t\ast}\L)= k([\L^\ast])$$ where
$[\L]\in\X_t$ is the isomorphism class of $\L$. By Riemann-Roch we get
$-\hat w = \f (\mu)$. (This implies $\f '(\hat\mu)=w$; after swapping $X$
and
$\X$ we get
$\f '(\hat\mu)=-w$ which implies $\f (w)=\hat\mu$.) 

To prove 4 we apply Riemann-Roch to
$$\bS ^0(\O_H)=\O_{\X},\qquad \bS ^1(\O_H)=0\,. $$ 5 is now
straightforward. Using these results one proves 6 as in \cite{Muk4}. 7
follows from 1 and 6.
\end{pf} If one defines a modified,
$H^\bullet(\cpuno,\Q)$-valued Mukai pairing by letting $\alpha\cdot \alpha'= 
p_\ast(\alpha^\ast
\cup\alpha')$ then the map $\f $ establishes an isometry
$H^\bullet(X,\Q)\iso H^\bullet(\X,\Q)$ as $H^\bullet(\cpuno,\Q)$-modules.
\begin{prop}
For all $\alpha\in H^\bullet(X,\Q)$, the $H^0(\X,\Q)$-component of $\f
(\alpha)$ is $\mu\cdot
\alpha$. As a consequence, $\f $ induces an isometry $\tilde \f \colon
\mu^\perp/\Q\mu\to H^2(\X,\Q)$. \label{ftilde}\end{prop} \begin{pf} We
already know that $\f (w)^0=0$ and $\f (1)^0=0$, so we may assume
$\alpha\in H^2(X,\Q)$. Then,
$\f (\alpha)^0=\pi^\ast\alpha\,\backslash\,\mu=\alpha\cdot\mu$.  Thus $\f
(\alpha)^0=0$ for
$\alpha\in\mu^\perp$. We now define $\bar \f \colon \mu^\perp\to
H^2(\X,\Q)$ by taking
$\bar \f (\alpha)$ as the $H^2$-component of $\f (\alpha)$. One has that
$\bar \f (\alpha)=0$ 
if and only if $\f (\alpha)=s\hat w$ ($s\in\Q$), and then
$\alpha=-s\hat\mu$, which proves
that $\ker\bar \f =\Q\mu$, and $\bar \f $ induces an injective morphism
$\tilde \f \colon
\mu^\perp/\Q\mu\hookrightarrow H^2(\X,\Q)$. If $\beta\in H^2(\X,\Q)$, $\f
'(\beta)\cdot \mu=0$, and $\beta=\tilde \f (-\f '(\beta))$, thus finishing
the proof.
\end{pf}
\begin{remark} The cohomology lattice $H^\bullet(X,\Z)$ contains a hyperbolic sublattice
$U$ generated by $\mu$ and $H$, and the hyperbolic sublattice
$V=H^0(X,\Z)\oplus H^4(X,\Z)$.
{From} Proposition \ref{ftilde} we see that (after identifying $X$ and
$\X$) the map $\f$ swaps the lattices $U$ and $V$. \end{remark} 

\smallskip {\it Topological invariants of the Fourier-Mukai transform.} Let
us assume at first that the Picard number of $X$ is two; then the Picard
group of $X$ reduces to the hyperbolic lattice $U$ (this happens when $X$
has 24 singular 
fibres consisting in elliptic curves with a nodal singularity). It is then
possible to compute
the invariants of the Fourier-Mukai transform of a sheaf on $X$ by means of
the Riemann-Roch
formula, expressed in the form $$\ch\bT ^\bullet(\F)=\frac1{\stdhp}\f
((\ch\F)\stdp)\,.$$ In
particular, let us assume that $\F$ is WIT$_i$, and set $\widehat
\F=\bT ^i(\F)$, and
$$\ch\F=r+a\,H+b\,\mu+c\,w\qquad\text{where}\quad r=\o{rk}\F\,.$$ We then have
\begin{equation} (-1)^i\,\o{rk}\widehat\F=a,\qquad
(-1)^i\,c_1(\widehat\F)=-r\,\Theta+c\,\hat\mu,\qquad (-1)^i\,
\ch_2(\widehat\F)=-b\,\hat
w\,.\label{e:RR1} \end{equation} In the same way, if
$\E$ is a WIT$_i$ sheaf on $\X$, with $$\ch\E=r+a\,\Theta+b\,\hat\mu+c\,
\hat w$$ after 
setting
$\widehat \E = \widehat{\bT}^i\E$ we have \begin{equation}
(-1)^i\,\o{rk}\widehat\E=a, \qquad
(-1)^i\,c_1(\widehat\E)=-r\,H+c\,\mu,\qquad (-1)^i\,
\ch_2(\widehat\E)=-b\,w\,.\label{e:RR2} \end{equation} One obtains similar
formulae also in the case when the Picard group has higher rank; in the
Appendix we treat the case when $X$ has also singular fibres of type
$I_n$ (according to Kodaira's classification \cite{Kod}). 

\oskip\section{Fourier-Mukai functor as mirror symmetry} We would like now
to examine some facts which pinpoint the relations between the relative
Fourier-Mukai transform on elliptic K3 surfaces and mirror symmetry.

\smallskip

(a) The formulae (\ref{e:RR1}) and (\ref{e:RR2}) establish a morphism
\begin{equation}
\begin{align*} H^0(X,\Z)\oplus \Pic(X) \oplus H^4(X,\Z) & \to H^0(\X,\Z)
\oplus \Pic(\X) \oplus H^4(\X,\Z) \\ r + a\,H + b\,\mu + c\,w & \mapsto a -
r\,\Theta + c\,\hat\mu -b\,
\hat w
\end{align*}\end{equation} together with its inverse. According to these
formulae, the cycle corresponding to a 0-brane is mapped to a special
Lagrangian 2-cycle of genus 1 (i.e.~to the cycle homologous to $\hat\mu$),
and {\it vice versa,} while a 4-brane is mapped to a special Lagrangian
2-cycle of genus 0, and {\it vice versa.} So one recovers the
transformation properties of D-branes under T-duality as known from string
theory
\cite{OOY}. One should notice that, according to Corollary \ref{cor2}, a
fibre of $X$, regarded as supersymmetric 2-cycle, is mapped to 0-brane
(point) lying in the same fibre, thus giving rise to a relative (fibrewise)
T-duality.

\smallskip

(b) Mirror symmetry should consist in the identification of the moduli
space of complex structures on an $n$-dimensional Calabi-Yau manifold $X$
with the moduli space of ``complexified K\"ahler structures'' on the mirror
manifold $\X$. The tangent spaces to the two moduli spaces are the
cohomology groups $H^{n-1,1}(X,\C)$ and $H^{1,1}(\X,\C)$, respectively. We
want to show that when $X$ is an (algebraic) elliptic K3 surface the $\f $
map establishes an isometry between the subspaces of these tangent spaces
which describe ``algebraic deformations,'' in a sense that we shall clarify
hereunder. 

We denote by
$\phi$ the complexification of $\tilde \f $ and by $\psi\colon H^2
(\X,\C)\to H^\bullet(X,\C)$ its inverse. 

\begin{prop} The map $\psi$ establishes an isometry $$\frac{H^{1,1}
(\X,\C)}{\o{Pic}(\X)\otimes\C}
\iso\frac{H^{1,1}(X,\C)}{\o{Pic}(X)\otimes\C}$$ \end{prop} \begin{pf} Let
$\Omega$, $\bar\Omega$ be generators of $H^{2,0}(\X,\C)$ and
$H^{0,2}(\X,\C)$. Since $\X$ is a moduli space of sheaves on $X$, by Remark
\ref{sympl} the classes
$\psi(\Omega)$ and
$\psi(\bar\Omega)$ lie in
$H^{2,0}(X,\C)$ and $H^{0,2}(X,\C)$, respectively. The result then follows
from Lemma
\ref{newlemma} and Proposition \ref{ftilde}. \end{pf} The space ${H^{1,1}(\X,\C)}/{\o{Pic}(\X)
\otimes\C}$ may be naturally identified with the tangent space at $\X$ to
the space of deformations of algebraic structures on $\X$ which preserve
the Picard lattice. Analogously, the space ${H^{1,1}(X,\C)}$ can be
regarded as the space of deformations of the K\"ahler structure of $X$, and
its quotient ${H^{1,1}(X,\C)}/{\o{Pic}(X)\otimes\C}$ as the space of
deformations of the K\"ahler structure which preserve the Picard lattice.
The map $\psi$ can then be thought of as a mirror transformation in the
algebraic setting. Since the Weil-Petersson metrics on both spaces are
expressed in terms of the Mukai pairing, which is preserved by $\psi$, we
see that $\psi$ establishes an isometry between the tangent spaces to the
two moduli spaces, consistently with the fact that the quantum cohomology
of a K3 surface is trivial.

\smallskip

(c) The mass of a BPS state, which is represented by a D-brane wrapped
around a 2-cycle $\gamma$, is given by the expression \cite{GK}
$$M=\frac{\left\vert\int_\gamma\Omega\right\vert}
{\left(\int_X\Omega\wedge\bar\Omega\right)^{\frac12}}=
\frac{\left\vert\gamma\cdot[\Omega]\right\vert}{\left([\Omega]
\cdot[\bar\Omega]
\right)^{\frac12}}$$
where $\Omega$ denotes a holomorphic 2-form on $X$, and $[\Omega]$ its
cohomology class in
$H^{2,0}(X,\C)$. The map $\phi$ evidently preserves this quantity. 

\smallskip

As a final remark, we would like to mention \cite{HO}, where the authors
consider a Fourier-Mukai transform on the cartesian product $X\times\X$
given by the ideal sheaf of the diagonal and use it to define a T-duality
between $X$ and
$\X$. A Riemann-Roch computation is then advocated to support an
interpretation of the duality of the baryonic phases in $N=2$ super
Yang-Mills theory. Thus, the geometric setting and the physical
implications of this construction are different from those of the present
paper.

\smallskip {\it Conclusions.} It should be stressed that in this picture,
in accordance
with
\cite{SYZ,Morr}, and differently to other proposals that have been recently
advocated (cf.~e.g.~\cite{Asp,GW,G}), the mirror dual to a given elliptic
K3 surface $X$ is isomorphic to $X$. Of course this does not imply that the
mirror map is to be trivial, and indeed the Fourier-Mukai transform seems
to establish such a map, at least at cohomological level, and at the
``infinitesimal'' level as far as the moduli spaces of complex structures
and the moduli space of complexified K\"ahler structures are concerned. It
would be now of some interest to develop a similar construction in terms of
a generalized Fourier-Mukai transform in higher dimensional cases, where
the mirror dual is not expected to be isomorphic to the original variety;
however, one may conjecture that the derived categories of the two
varieties are equivalent.

\medskip\noindent {\bf Acknowledgements.} We thank C.~G\'omez, C.-S.~Chu,
and especially C.~Imbimbo for useful discussions. This research was partly
supported by the Spanish DGES through the research project PB95-0928, by
the Italian Ministry for Universities and Research, and by an
Italian-Spanish cooperation project. The first author thanks the Tata
Institute for Fundamental Research, Bombay, for the very warm hospitality
and for providing support during the final stage of preparation of this
paper.

\oskip\section{Appendix} In order to be able to compute the topological
invariants of the Fourier-Mukai transform of a sheaf on $X$ we need to
describe the action of the $\f $ map on the generators of the Picard group.
In this Appendix we assume that the elliptic K3 surface $X$ has singular
fibres which are of type $I_n$, $n \ge 3$ or are elliptic nodal curves;
every singular fibre of type $I_n$ is a reducible curve whose irreducible
components are
$n$ smooth rational curves which intersect pairwise. The section $e$
intersects only one irreducible component of each singular fibre. The
Picard group $\pic X$ is generated by the divisors $\mu$ and $H$ and by $r$
divisors
$\alpha_1,\dots,\alpha_r$ given by the irreducible components $C_i$ of the singular fibres of type $I_n$ which do not meet the section $e$. 

Since $\alpha_i\cdot\mu=0$ and $\alpha_i\cdot H = \f (\alpha_i)\cdot(1+\hat
w)=0$ we have
$\f (\alpha_i)=\beta_i\in H^2(\X,\Q)$.

\begin{prop} The sheaf $\O_X(-C_i)$ is WIT$_1$, and $\bT^1
\left[\O_X(-C_i)\right]\simeq
\O_{\Sigma_i}(-1)$, where $\Sigma_i$ is a section of $\X$ whose associated
cohomology class is $\Theta+\hat\mu+\beta_i$. \end{prop} \begin{pf} By base
change for every $t\in\cpuno$ one has \begin{equation}
\bT^i\left[\O_X(-C_i)\right]\otimes\O_{\X_t}\simeq
\bT^i_t(\L_t)\label{e:rf}
\end{equation} where $\L_t=\rest{\O_X(-C_i)},{X_t}$ so that $\O_X(-C_i)$ is
WIT$_1$. By Riemann-Roch one has
\begin{equation}\ch\bT^1\left[\O_X(-C_i)\right]=\Theta+\hat\mu+\beta_i
\label{e:LRR}\end{equation} {From} equation (\ref{e:rf}) we see that $\bT^i
\left[\O_X(-C_i)\right]
\otimes\O_{\X_t}$ is concentrated at the point in $\X_t$ corresponding to
the flat line bundle $\L_t^\ast$ on $X_t$, whence the first claim follows.
The second is a consequence of formula (\ref{e:LRR}). \end{pf}

\oskip

\end{document}